\documentclass[review]{elsarticle}

\usepackage{lineno,hyperref}
\modulolinenumbers[5]

\journal{Communications in Mathematical Physics}

\begin{document}

\begin{frontmatter}

\title{Physical Versus Mathematical Billiards: From Regular Dynamics to Chaos and Back}



\author{L.A.Bunimovich\corref{mycorrespondingauthor}}
\cortext[mycorrespondingauthor]{Corresponding author}
\ead{leonid.bunimovich@math.gatech.edu}

\address{School of Mathematics,
Georgia Institute of Technology,
Atlanta, GA 30332-0160 USA}

\begin{abstract}
In standard (mathematical) billiards a point particle moves uniformly in a billiard table with elastic reflections off the boundary. We show that in transition from mathematical billiards to physical billiards, where a finite size hard sphere moves in the same billiard table, virtually anything may happen. Namely a non-chaotic billiard may become chaotic and vice versa. Moreover, both these transitions may occur softly, i.e. for any (arbitrarily small) positive value of the radius of a physical particle, as well as by a "hard" transition when radius of the physical particle must exceed some critical strictly positive value. Such transitions may change a phase portrait of a mathematical billiard locally as well as completely (globally). These results are somewhat unexpected because for all standard examples of billiards their dynamics remains absolutely the same after transition from a point particle to a finite size ("physical") particle. Moreover we show that a character of dynamics may change several times when the size of the particle is increasing.
\end{abstract}


\end{frontmatter}

\section*{Introduction}

Billiards are dynamical systems generated by an uniform motion of a point particle within a domain with a piecewise smooth boundary.
Upon reaching the boundary of the domain (billiard table) the particle gets elastically reflected. 
It is a commonly held opinion that such mathematical billiards generated by the motion of a point particle adequately describe dynamics of real physical particles within the same billiard table. By a physical particle we mean here a hard sphere of radius $r$ which gets elastically reflected off the boundary of a billiard table. 
Clearly one can follow evolution of a (spherical) physical particle by considering motion of its center. Therefore it is easy to see that dynamics of a physical billiard is equivalent to dynamics of mathematical billiard within a smaller billiard table which is obtained by shrinking the initial billiard table by moving points of its boundary on the distance $r$ along normal vectors pointing towards the interior of the table. (Clearly a particle cannot move within a billiard table at all if it is too big for this table. Therefore we always assume that radius of a moving sphere allows it to move within some sub-domain of the billiard table, which becomes a billiard table for the mathematical billiard  corresponding to our physical billiard).
One can immediately see that for all basic and most popular classes of billiards transition from mathematical billiards to physical billiards makes absolutely no change to dynamics (besides, of course, making all free passes shorter). Indeed a physical billiard within a triangle becomes mathematical billiard within a smaller similar triangle. Analogously a physical billiard within a circle becomes a mathematical billiard within a smaller circle, Sinai billiard remains a Sinai billiard with larger scatterers and squash billiards result in smaller squashes. (Recall that the boundary of squash billiards,(sometimes called tilted stadia), consists of two circular arcs connected by two tangent to them straight segments. A squash becomes a stadium if the arcs are semicircles of the same radius. 
This is a reason for a general opinion that transition to physical billiards is not interesting because it brings nothing new. In fact, when some new question(s) arise we always look what is going to happen in the most popular and visual examples.
In this paper we show that this general opinion is completely wrong. Moreover we demonstrate that in fact virtually any possible changes of dynamics may occur as a result of transition to a physical billiard. Namely such transition may generate appearance of KAM-tori when mathematical billiard was completely chaotic. The opposite transition, when a non-chaotic mathematical billiard becomes completely chaotic physical billiard may also happen.
Moreover, both such transitions may occur for any positive $r>0$, i.e. demonstrating a "soft" transition, as well as, via 'hard" transitions when a value of the radius of the moving physical particle must exceed some strictly positive critical value.
We present concrete examples of all types of transitions. 
All the examples we consider are two-dimensional because a goal of this paper is just a proof of a concept. Therefore we choose and build the most visual (and hopefully the simplest) examples. 
This paper opens up a new direction in the studies of billiards. In this connection it is worth to remind readers who doubt that because, e.g. hard spheres gases could be reduced to mathematical billiards, that we consider here motion of just one particle (sphere). In fact in order to analyze dynamics of physical billiards we always reduce them to the mathematical ones. However in a course of such reduction the shape of a billiard table may change essentially.
This paper opens up a new area not just for mathematical but for physical studies as well. For instance, real pipes, containers and channels always have non-ideally smooth boundaries and generally non piecewise  smooth boundaries as well. In fact the "real" boundaries are rough and contain impurities of different scales. Our examples show that generally such impurities may slow down a flow by generating local vortices because of appearance of KAM islands which lead to destruction of regions with chaotic dynamics. Analogously some islands (vortices) may disappear in transitions to physical billiards. This observation should be taken into consideration both in theoretical and numerical studies, especially when considering propagation in nanochannels \cite{Ro} where a width of a channel is very small and therefore particles propagate one after another (i.e. do not change their "order"). 
Especially the studies of physical billiards are relevant for quantum mechanics because of the Heisenberg's uncertainty principle. Hence (mathematical) point particles do not seem to be relevant in the studies of quantum chaos. Instead one should consider evolution of wave packets, i.e. of finite size "particles". Therefore a quantum billiard may be more chaotic than the corresponding classical billiard. Remarkably, this claim made by the author at a recent conference inspired physicists, and, it was just shown \cite{RBG} that an uniformly held opinion that classical systems are always "more chaotic" than their quantum counterparts is not always true. An opposite situation may occur as well. 
\section{Physical Billiards}
Consider a free motion of a (hard) sphere with radius $r$ within a domain $Q$ in $d$-dimensional Euclidean space. We always assume that the boundary $\partial Q$ of $Q$ consists of a finite number of co-dimension one smooth manifolds of class at least $C^2$ which are called regular components of the boundary. Interior points of regular components are called regular points of the boundary, and points of intersection of regular components will be called corner (or singular) points. In what follows such domain $Q$ will be often called a billiard table. 

At any interior point $q$ of a regular component of the boundary $\partial Q$ there exists an unique inner (i.e. toward the interior of $Q$) unit normal vector $n(q)$. Upon reaching the boundary of a billiard table the particle gets elastically reflected off it, i.e. the angle between the velocity vector $v$ of the particle and $n(q)$ before reflection (the angle of incidence) equals the angle between $n(q)$ and the velocity of the particle after reflection (the angle of reflection). 

One gets a standard (mathematical) billiard if a moving particle is a point, i.e. its "radius" $r=0$. 
We call a billiard physical if a moving particle is a hard sphere with a positive radius.
Physical (as well as mathematical) billiards are Hamiltonian systems. Therefore these dynamical systems have a natural invariant measure which is a volume in the phase space.
We will consider a billiard map which arises if one follows billiard trajectories only at the moments immediately after the particle's reflections off the boundary of a billiard table.
Such billiard map generates a 
dynamical system (see e.g. \cite{CM}). For our purposes in this paper neither an exact expression for a billiard map nor coordinates in a corresponding phase space are needed. Therefore we will stick with just visual purely geometric consideration.  
A regular component of the boundary $\partial Q$ is called dispersing (focusing) if it is convex inwards (outwards) of a billiard table (Fig.1). We assume in what follows that a curvature of dispersing (focusing) components is negative (positive) in all their points. A regular boundary component is called neutral if its curvature is identically zero.
It is easy to see that the moving particle with a positive radius $r$ may never hit some parts of the boundary of the billiard table. For instance in case when dimension of the billiard table is two it happens if at least two regular components of the boundary intersect under the angle less than $\pi$. Indeed in such case a disk (particle) with positive radius can not get into this corner and hit points at the boundary of the billiard table which are close to a point of intersection of these regular components.
\begin{figure}[ht]
    \centering
    \includegraphics[scale=0.25]{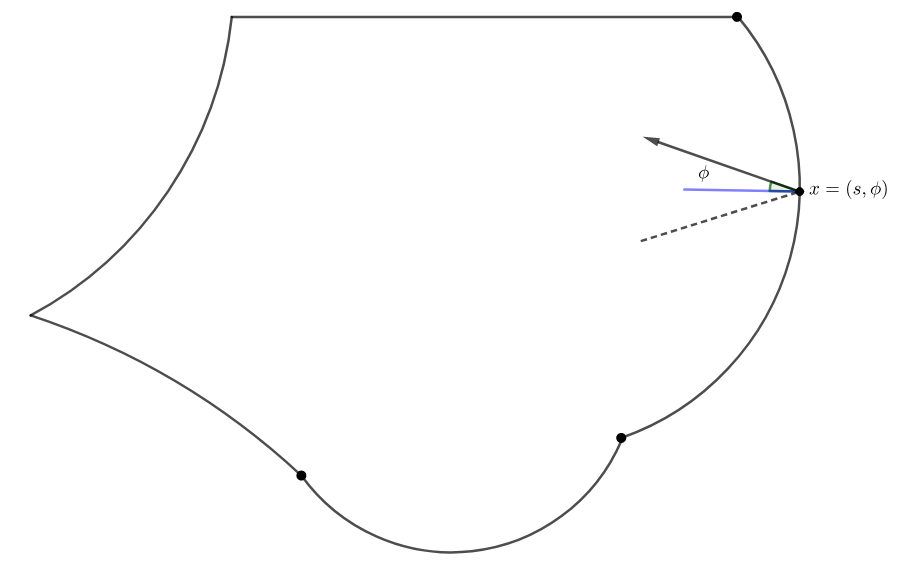}
    \caption{A billiard table.}
\end{figure}
We now formally describe process of transition from a mathematical billiard to the corresponding physical one.  
Consider a physical (finite size) particle moves in the same billiard table as a point particle in the mathematical billiard, to which we will compare the physical one. 
To represent dynamics of a (hard) homogeneous spherical particle of radius $r$ it is enough to follow motion of its center. It is easy to see that the center of particle moves in the smaller table which one gets by moving all points $q$ of the boundary by $r$ into the interior of the billiard table along the internal normal vector $n(q)$ (see Fig.2). Dynamics of the center of a physical particle is equivalent to a mathematical billiard in this smaller billiard table which we will call from now on a reduced billiard table.  
We assume throughout this paper that billiard tables satisfy the following Condition NIC, where NIC is an abbreviation for "no internal corners".

Definition. A billiard table $Q$ satisfies Condition NIC if either 

(i)all straight segments which connect any point on any regular component $\gamma$ of the boundary $\partial Q$ to any point of any regular component intersecting $\gamma$ belong to $Q$

or, if (i) does not hold, then

(ii)any two intersecting regular components of the boundary $\partial Q$ for which condition (i) does not hold have a common tangent at their point of intersection. 
In other words the boundary $\partial Q$ is of class $C^1$ at such points of intersection of regular components. 
Here NIC is an abbreviation for "no internal corners".
\begin{figure}[ht]
    \centering
    \includegraphics[scale=0.25]{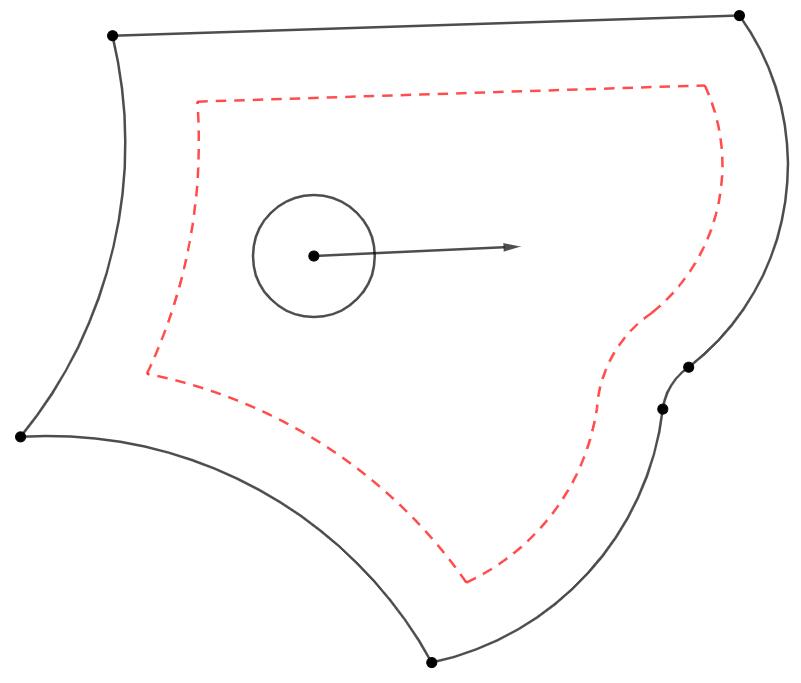}
    \caption{Transition from mathematical to a physical billiard.}
\end{figure}
For example the billiard table in Fig.1 does not satisfy Condition NIC because some segments connecting points on two adjacent focusing components do not belong to the billiard table. On another hand the billiard table in Fig.2 satisfies condition NIC. Observe also that this condition does not allow billiard tables to be non-convex polygons.
This restriction (condition NIC) will be lifted in a consequent paper, but here we do not want to deal with technical complications which may confuse otherwise clear and visual transition from physical billiards to dynamically equivalent mathematical billiards.
For the same reason we will consider in what follows only two dimensional billiard tables.
It is possible to give many other examples demonstrating our main claim that in transitions from mathematical to physical billiards various changes in dynamics may occur. However a purpose in this paper is just a proof of a concept. Therefore our examples will be the most simple and visual ones. Moreover in order to avoid long technical proofs typical for the theory of billiards all statements about changes in dynamics which occur in transitions from mathematical to physical billiards  
follow from known results about mathematical billiards.

\section{Birth and destruction of KAM-islands in transition from mathematical to physical billiards}
In this section we consider appearance and disappearance of KAM-islands. These examples will demonstrate that transitions from mathematical to physical billiards may result in changes of the phase portrait of the corresponding dynamical system.
Following our general strategy in this paper only the simplest situation will be discussed when KAM-islands are generated by elliptic periodic orbits of period two. Indeed, it is the simplest case because in billiards there are no fixed points.
Necessary and sufficient conditions for linear stability of period two points in billiards are well known \cite{W}. They read as

\begin{equation}
L k_0 k_1-k_0-k_1<0
\end{equation}

\begin{equation}
(Lk_0-1)(Lk_1-1)>0    
\end{equation}

where $k_1$ and $k_2$ are curvatures of the boundary at the endpoints of a period two orbit, and $L$ is the length of the segment connecting these points.
It can be immediately seen that if the curvatures $k_0, k_1$ have the same sign then the condition (1) is violated for large $L$, and if $k_0, k_1$ have opposite signs then the condition (2) is violated for large $L$. Observe also that an orbit ending on two dispersing, or on two neutral components is always unstable. 

The case when only one end of orbit belongs to a neutral component can be reduced to consideration of orbit bouncing between two focusing or two dispersing components by a standard trick of reflecting a billiard table with respect to a neutral component of the boundary.

Therefore only the cases when both corresponding components are focusing or one is focusing and another one is dispersing are of interest.

Consider first period two orbits bouncing between two focusing components. Let $q_1$ and $q_2$ are the ends of a periodic orbit and regular points of the boundary components of some billiard table. Denote by $R_1$($R_2$) a radius of curvature of the first (second) focusing component at the point $q_1$($q_2$). Without loss of generality we assume that $R_1>R_2$. Then the conditions of linear stability (1) and (2) become $L<R_1+R_2$ and $(L-R_1)(L-R_2)>0$ respectively.
Suppose that this periodic orbit is linearly stable and $L<R_1+R_2$, $L>R_1$, $L>R_2$. 

Let a physical particle with radius $r$ moves in a corresponding billiard table. Then in the reduced (by $r$ inside) billiard table this period two orbit remains. However, its length becomes $L-2r$ while both radia of curvature at its end points will become $R_1-r$ and $R_2-r$ respectively.
Therefore this orbit becomes unstable when $L<R_1-r$ but $L>R_2-r$. In other words period two orbit loses stability when radius $r$ of a physical particle equals $L-R_1$. If radius is smaller than this value then this orbit remains linearly stable.

However if radius of the particle continues to increase then after passing the value ($L-R_2$) another change in dynamics (bifurcation) occurs. Indeed it follows from (2) that period two orbit acquires linear stability if the radius of moving particle exceeds $L-R_2$. (Of course $r$ must remain to be less than $L/2$, and thus $L$ must be less than $2R_2$, otherwise this orbit will just disappear).

Therefore by increasing the radius of moving particle it is possible to achieve several changes in dynamics (or, more precisely, in the phase portrait) of the corresponding billiard. Indeed, in our example a periodic orbit is stable if radius of moving particle is relatively small, then it becomes unstable when this radius is within a range of larger values. And finally an orbit again acquires stability when the radius of the particle becomes large enough.
In fact a radius $r$ of the physical particle can be viewed as 
a bifurcation parameter in this family of dynamical systems (physical billiards).

However, as we will see in the next section, variations of the radius of physical particle may lead to enormous changes of the corresponding billiard tables which are not really bifurcations.
It is well known \cite{Mo} that linearly stable periodic orbits not always generate KAM-islands. To make this happen some quantity called the first Birkhoff coefficient \cite{Mo} must not vanish. Again for simplicity we assume that both focusing components which we consider are arcs of some circles. In this case the condition of ellipticity is reduced to a few inequalities \cite{SdK1,SdK2}. Namely $4(L-R_1)(L-R_2)$ must not be equal to $R_1R_2$ or to $2R_1R_2$. It is very easy to satisfy these conditions by choosing appropriate values of parameters. 

Generally, it is well known that ellipticity is an open property, in the sense that sufficiently small $C^2$ perturbations of the focusing components near the ends of elliptic periodic orbit in billiards will again have an elliptic periodic orbit (close to the perturbed one).               
Thus just pick such radia of circles that the above two inequalities (conditions of ellipticity) hold. 

In the next example we will consider another type of appearance of a linearly  stable period two orbit from a linearly unstable one which occurs in transition from mathematical to physical billiards. 
Take now a billiard orbit which moves between focusing and dispersing components. We will keep the same notations as above. Now $R_1$ ($R_2$) is a radius of curvature of the focusing (dispersing) component of the boundary at the corresponding ends of an 
orbit.
In this case conditions of stability (1) and (2) become $L>R_1-R_2$ and $(L-R_1)(L+R_2)<0$. Therefore $L$ must be smaller than $R_1$. Suppose that in a mathematical billiard this orbit is linearly unstable because $L>R_1$. Then this orbit will become stable when radius $r$ of the particle exceeds $L-R_1$, and generically a new KAM-island will appear if $L<2R_1$ 

\section{Global transitions between non-chaotic and chaotic billiards when radius of the particle changes}

In the previous section we showed that phase portraits of mathematical billiards may change in the transition from mathematical to physical billiards. These changes resulted in appearance and disappearance of KAM-islands, i.e. in a sense, local ones. 
In this section we give examples of global transitions from non chaotic billiards to the fully chaotic ones and vice versa.  By fully chaotic dynamical systems we mean the ones which are completely hyperbolic. Recall that a system is completely hyperbolic if through almost any point in its phase space pass smooth (local) stable and unstable manifolds. Thus, in other words by fully chaotic systems are dynamical systems with nonvanishing Lyapunov exponents. 
It is well known that such systems have a finite or countable number of ergodic components of positive measure which cover (up to a set of measure zero) the entire phase space. Moreover, these dynamical systems have positive Kolmogorov-Sinai entropy, are mixing on each ergodic component and have as well other strong chaotic properties (see e.g. \cite{S,KH}).

For the sake of clarity we will present only simple and visual examples of such transitions. Moreover proofs of the corresponding claims will directly follow from already known (rigorous) results on mathematical billiards. 
\begin{figure}[ht]
    \centering
    \includegraphics[scale=0.25]{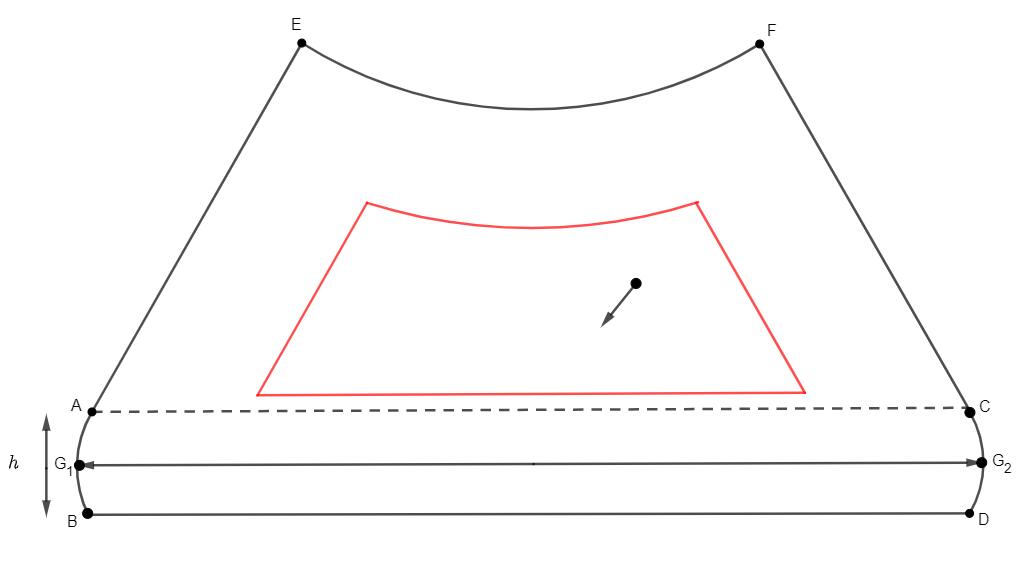}
    \caption{A hard transition to chaotic billiard.}
\end{figure}

Consider first a billiard table $Q$ depicted in Fig.3. Its boundary consists of six regular components. Two of them are focusing and formed by equal arcs AB and CD of some circles of the same radius $R$. The neutral components AE, BD and CF are straight segments tangent to both these arcs at the points A and C. The last smooth component EF of the boundary $\partial Q$ is dispersing. We also assume that the angles between the segments AE and AC and between the segments AC and CF are both equal $\pi/3$ and that the segments AC and BD are parallel.
Besides the arcs AC and BD are placed so that their centers $G_1$ and $G_2$ are at a distance strictly less than $2R$ from each other and the segment $G_1G_2$ contains the centers of both circles to which the arcs AB and CD belong. Then the period two orbit $G_1G_2$ is linearly stable. Conditions of nonlinear stability of such orbit are given in the previous section. We assume that these conditions are satisfied and thus our orbit is the center of a KAM-island. It is also assumes that dispersing component EF is placed at such distance from the segment BD which is larger than twice the distance between parallel segments AC and BD.
Let now a hard disk with radius $r$ moves within $Q$ and generates a physical billiard.

Lemma 1. A physical billiard within $Q$ is completely hyperbolic (chaotic) if the radius of the moving particle exceeds the length of the straight segment connecting points A and B.  

Proof. Let the radius of the moving particle is larger that the segment AB (which, according to our conditions, equals the segment CD). Then physical billiard within the billiard table ABDCFE is equivalent to mathematical billiard with a smaller billiard table depicted in Fig.3. 
The billiard table of this smaller mathematical billiard contains three neutral components and one dispersing component of the boundary. Indeed under the conditions of the lemma elliptic island completely disappears when radius $r$ of the particle is greater than the length of segment AB. (In fact this elliptic island disappears when $r$ becomes greater than half of the length of the segment AB. However, then some extra computations are needed to show that the corresponding mathematical billiard is hyperbolic. We want the results though to follow from already known facts about mathematical billiards in order not to corrupt clear ideas by extra computations and formulas. Proof of a concept is a key goal.)  

Observe now that in view of the conditions on the initial large billiard table ABDCFE three neutral components of the small billiard table in Fig.3 belong to a regular triangle. The last forth regular component of the boundary of the small table is dispersing. Therefore this mathematical billiard dynamically equivalent to physical billiard if radius of the particle is greater than the length of the segment AB. It is well known (see e.g. \cite{CM} and references therein) that this mathematical billiard is completely hyperbolic.  
Clearly an elliptic orbit in the phase space of a mathematical billiard generated by our physical billiard remains for all $r<|AC|/2$. 
Therefore in this example a birth of completely chaotic billiard in the transition from mathematical billiard to physical billiard occurs in a "hard way", i.e. when radius of physical particle exceeds some critical value. 
\begin{figure}[ht]
    \centering
    \includegraphics[scale=0.25]{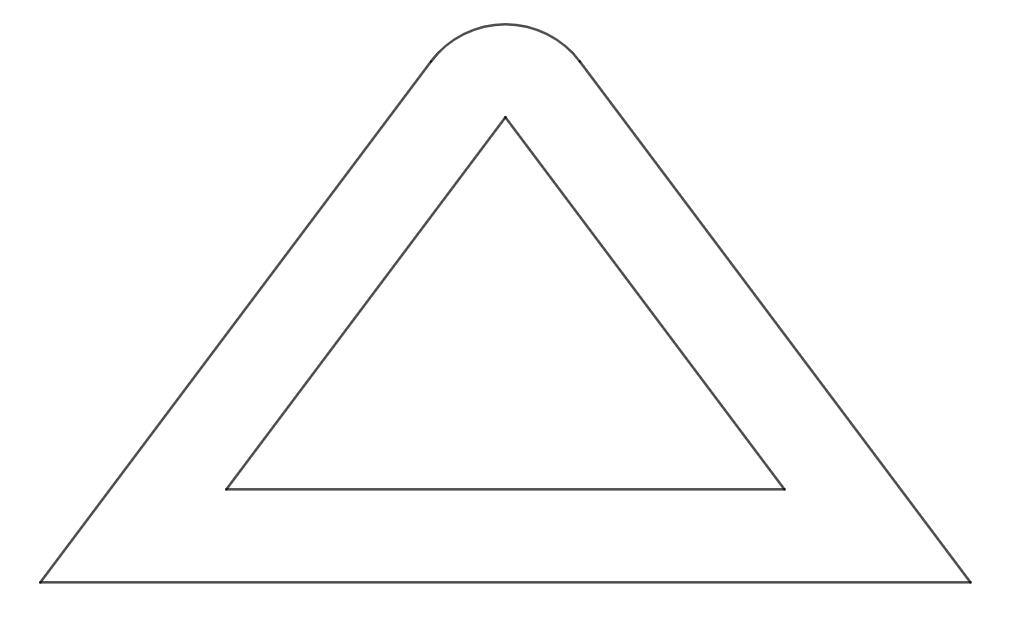}
    \caption{A hard transition to non-chaotic billiard.}
\end{figure}

We now consider opposite transition from (completely) chaotic mathematical billiard to non-chaotic physical billiard. Moreover, the resulting physical billiard will be strongly non-chaotic. More precisely, this physical billiard does not have any subset with hyperbolic dynamics and its Kolmogorov-Sinai entropy equals zero. 
Take a billiard table having a shape of a regular triangle with one vertex smoothened by an arc of a circle with radius $R$ (Fig.4). It is well known (see e.g. \cite{B1,B2}) that this mathematical billiard is strongly chaotic, i.e. completely hyperbolic. In fact, this biliard is also ergodic.
It is easy to see that physical billiard in this table becomes equivalent to a mathematical billiard in a regular triangle when radius of the particle becomes equal $R$.
Therefore we have

Lemma 2. There exist billiard tables with strongly chaotic mathematical billiards where physical billiards become completely non-chaotic if radius of the particle exceeds some critical value.
\begin{figure}[ht]
    \centering
    \includegraphics[scale=0.25]
    {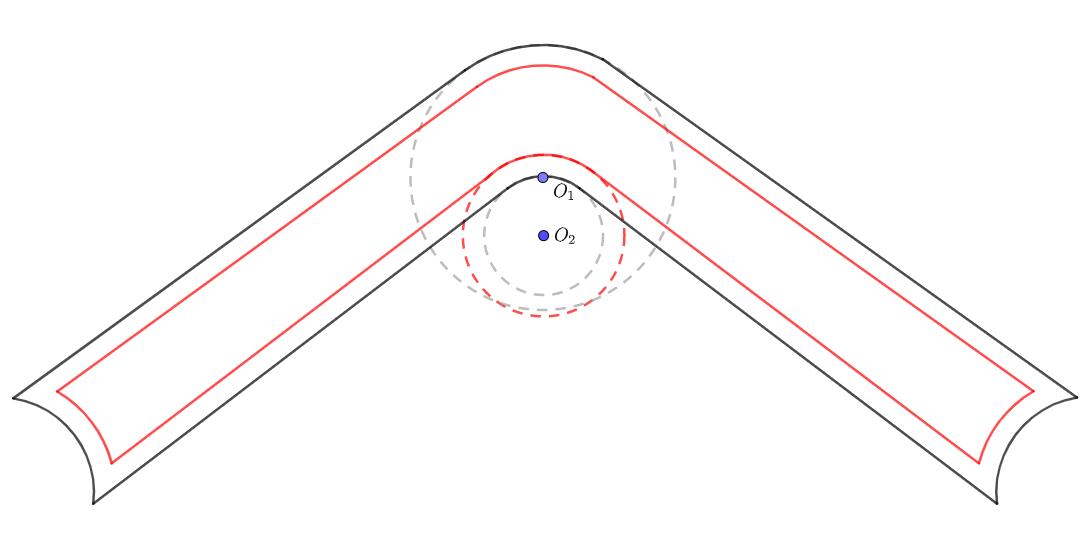}
    \caption{A soft transition to a non-chaotic billiard.}
\end{figure}

Now we give an example of another transition when physical billiards within billiard tables of chaotic mathematical billiards acquire a linearly stable periodic orbit (i.e. generically a KAM-island) and thus become non-chaotic. 
Consider a billiard table $Q$ depicted in Fig.5. The boundary $\partial Q$ of this billiard table contains eight regular components. Two regular components are arcs of the circles with centers $O_1$ and $O_2$. One of these components is focusing and another one is dispersing. Four components adjacent to these two are neutral ones tangent to focusing or dispersing component at their endpoints. We also assume that the circle with the center $O_1$ is tangent to the circle with the center $O_2$ (Fig.5). The last two regular components of the boundary are dispersing ones.
It follows from the paper on track billiards \cite{BdM} that this (mathematical) billiard is hyperbolic.  

Lemma 3. A physical billiard in $Q$ has a linearly stable periodic point for any radius $r>0$ of the moving particle.

Proof. Consider period two orbit which lies on the line containing the centers $O_1$ and $O_2$ of the corresponding circles. It follows from the previous section that this periodic orbit becomes linearly stable for a physical billiard with $r>0$. 
Indeed in the transition to physical billiard dispersing arc will remain dispersing but radius of the circle which contains it will increase by $r>0$ and therefore this circle will contain $O_1$. On another hand the focusing component becomes an arc of a smaller circle with the same center $O_1$.  
Therefore conditions (1), (2) imply that this orbit is linearly stable.
To the best of our knowledge exact conditions of ellipticity of such period two orbit are not computed explicitly as it was done for orbits bouncing between two focusing components \cite{SdK1,SdK2}. However an ellipticity (nonvanishing of Birkhoff coefficient) is an open condition. Therefore there exist such values of the radia of circles with the centers $O_1$ and $O_2$ that the corresponding  orbit is elliptic.
In this examples we see a soft transition which occurs for any positive radius $r$ of the physical particle. In this transition a KAM-island appears being  generated by period two stable orbit. 

It is natural to conjecture though that not just one island appears here but infinitely many KAM-islands which coexist with positive measure sets carrying a chaotic (hyperbolic) dynamics. Indeed, it is a virtually universal belief that generic Hamiltonian systems have such (often called divided) phase space. However this claim is not proved yet. 

 Our last example will demonstrate that transition to completely hyperbolic (chaotic) billiard may occur without drastic changes as "erasing" some boundary components of the mathematical billiard table in the two previous examples.
Moreover the next example demonstrates that transition from nonchaotic to chaotic billiard may occur even in billiards with convex tables.
Consider a convex billiard table bounded by arcs of two circles. We assume that the arc with smaller radius $R_2$ is larger than a half of the corresponding circle. Such billiard tables are called skewed lemons \cite{BZZ}.
Assume that an arc with the larger radius $R_1$ is such that the center of this larger circle lies within billiard table. Then the orbit which bouncing between centers of these focusing components is stable. 
In the previous section we saw that this orbit loses stability when radius of the moving particle exceeds some critical value. It was proved in \cite{BZZ, JZ} that skewed lemon billiards are completely hyperbolic if $R_1$ is sufficiently large.  
Once again, for a proof of concept we will make a visual trick.
Consider a chaotic skewed lemon (mathematical) billiard. Such billiards exist \cite{BZZ, JZ}. Assume now that the corresponding billiard table $Q$ is in fact a reduced table which appeared in transition from mathematical to a physical billiard with moving particle of radius $r$. Take now $r$ greater than the distance from the center of the circle with the larger radius to the billiard table $Q$. Then this center belongs to the billiard table of the original mathematical billiard. Therefore period two orbit in this original mathematical billiard is stable. However a physical billiard is completely hyperbolic (chaotic). Therefore in this example a mathematical nonchaotic billiard becomes chaotic when the radius of a physical particle exceeds a critical value.    
 
\section{Concluding remarks}
We demonstrated that various changes in dynamics may occur because of transition from mathematical to physical billiards. Moreover, it can be several such transitions when radius of the particle is changing. Then there exist several different ranges of sizes of the particle where dynamics of the corresponding billiard changes when radius of the particle crosses the boundaries between these ranges. 
Besides raising many questions for the future mathematical studies these results should be of interest to physics. Indeed a quantum "particle" has a finite size thanks to Heisenberg's uncertainty principle. 

Indeed, following the ideas developed in our paper, it was demonstrated that, in sharp contrast to conventional wisdom, quantum system can be more chaotic than its classical counterpart \cite{RBG}. So far it was always thought that only the opposite is possible, i.e. classical systems are always more chaotic than their quantum counterparts.

We also claim that propagation of particles in channels may have different properties depending on their size. Indeed a birth and destruction of KAM-islands results in appearance and disappearance of "vortices" in a flow respectively. As it was demonstrated in this paper, even several changes in dynamics may occur when a size of the moving particle varies. These phenomena should be especially of interest for flows in nanochannels (see e.g. \cite{Ro}).

\section{Acknowledgements}
This work was partially supported by the NSF grant DMS-1600568.


\begin{thebibliography}{99}

\bibitem{B1}
Bunimovich, L.A.: On ergodic properties of some billiards. Funct. Anal. and  Appl-s  8, 254-255 (1974)

\bibitem{B2}
Bunimovich, L.A.: On the Ergodic Properties of Nowhere Dispersing Billiards. Commun. Math. Phys. 65, 295-312 (1979)


\bibitem{BdM}
Bunimovich, L.A., Del Magno, G.: Track billiards. Commun. Math. Phys. 288, 699–713 (2009)

\bibitem{BZZ} 
Bunimovich,L.A., Zhang, H.K., Zhang, P.: On another edge of defocusing: hyperbolicity of skewed lemon billiard. Commun. Math. Phys. 341, 781-803 (2016) 

\bibitem{SdK1}
Carneiro,M.J.D., Kamphorst,S., Pinto De Carvalho, S.: Elliptic islands in strictly convex billiards. Erg. Th. and Dyn. Syst. 23 (3), 799-812 (2003)

\bibitem{CM}
Chernov, N., Markarian, R.: Chaotic Billiards. AMS Publ., RI (2006)

\bibitem{Ro}
Jepps, O.G., Rondoni,L.: Thermodynamics and complexity of simple transport phenomena.  J. Phys. A: Math. Gen. 39, 1311-1323 (2006)

\bibitem{JZ}
Jin, X., Zhang, P.: Hyperbolicity of asymmetric lemon billiards. arXiv:1902.08130 [math.DS]

\bibitem{SdK2} 
Kamphorst,S., Pinto de Carvalho,S.: The first Birkhoff coefficient and the stability of 2-periodic orbits of billiards. Experimental Math. 14, 299-306 (2005)

\bibitem{KH}
Katok,A.B., Hassellblatt, B.: Introduction to the modern theory of dynamical systems, Cambridge Univ. Press, Cambridge (1995)

\bibitem{Mo}
Moser, J.: Stable and random motions in dynamical systems: with special emphasis on celestial mechanics. Princeton Univ. Press, Princeton (1973)

\bibitem{RBG}
Rosenbaum, E.B., Bunimovich, L.A., Galitski, V.: Quantum chaos in classically non-chaotic systems. arXiv:1902.05466 [quant-ph]

\bibitem{S} 
Ya. G. Sinai, Introduction to ergodic theory. Princeton Univ. Press, Princeton (1976)


\bibitem{W}
Wojtkowski, M.: Principles for the design of billiards with nonvanishing Lyapunov exponents.  Commun. Math. Phys. 105, 391-414 (1986)
 

\end{thebibliography}

\end{document}